\input amstex

\define\sgn{\operatorname{sgn}} 
\define\p. v.{\operatorname{p. v.}} 
\documentstyle{amsppt}
\NoBlackBoxes 
\topmatter
\title 
Estimates for singular integrals and extrapolation     
\endtitle
\topmatter
\author
Shuichi Sato
\endauthor
\thanks 2000 Mathematics Subject Classification. Primary 42B20, 42B25 
\endthanks  
\abstract
In this note, we study singular integrals with 
rough kernels, which belong to a class of singular Radon transforms.  
We prove certain estimates for the singular integrals that 
are useful in an  extrapolation argument.    As an application, 
we prove $L^p$ boundedness of the singular integrals under 
a certain sharp size condition on their kernels.  
\endabstract
\endtopmatter
\document
\head 1.  Introduction \endhead 
Let $\Omega$ be a function in $L^1(S^{n-1})$ satisfying 
$$ \int_{S^{n-1}} \Omega(\theta)\, d\sigma(\theta)=0, \tag 1.1 $$  
where $d\sigma$ denotes the Lebesgue surface measure on the unit sphere $S^{n-1}$ 
in $\Bbb R^n$. In this note we assume $n\geq 2$.  
 For $s\geq 1$, let $\Delta_s$ denote the collection of measurable functions 
$h$ on $\Bbb R_+=\{t\in \Bbb R : t> 0\}$ satisfying 
$$\|h\|_{\Delta_s}=  
\sup_{j\in \Bbb Z}\left(\int_{2^j}^{2^{j+1}}|h(t)|^s\,dt/t\right)^{1/s}
<\infty,$$ 
where $\Bbb Z$ denotes the set of integers.   
We note that $\Delta_s\subset \Delta_t$ if $s>t$.  In this note 
we always assume $h\in \Delta_1$.    
Let $P(y)=(P_1(y),P_2(y),\dots,P_d(y))$ be a polynomial 
mapping, where each $P_j$ is a real-valued polynomial on $\Bbb R^n$.   
We consider a singular integral operator of the form:  
$$ T(f)(x)=\p. v. \int_{\Bbb R^n}f(x-P(y))K(y)\,dy=
\lim_{\epsilon\to 0}\int_{|y|>\epsilon}f(x-P(y))K(y)\,dy,  \tag 1.2 $$  
for an appropriate function $f$ on $\Bbb R^d$, 
where $K(y)=h(|y|)\Omega(y')|y|^{-n}$,  $y'=|y|^{-1}y$.  Then,   $T(f)$  
 belongs to a class of singular Radon transforms. 
See Stein \cite{17}, Fan-Pan \cite{8} and Al-Salman-Pan \cite{1} for this 
singular integral.   
\par 
When $h=1$ (a constant function), $n=d$ and $P(y)=y$, we also 
write $T(f)=S(f)$.   
Let 
$\hat{f}(\xi)=\int_{\Bbb R^d} f(x)e^{-2\pi i\langle x, \xi\rangle}\,dx $
be the Fourier transform of $f$, where $\langle \cdot,\cdot\rangle$ denotes the inner product in $\Bbb R^d$.  
Then it is known that $(Sf)\hat{\phantom{t}}(\xi)=m(\xi')\hat{f}(\xi)$, 
where  
$$m(\xi')=-\int_{S^{n-1}}\Omega(\theta)\left[i\frac{\pi}{2}\sgn(\langle 
\xi',\theta\rangle)+\log |\langle \xi',\theta\rangle|\right]\, d\sigma(\theta).
$$  
Using this,  we can show that 
 $S$  extends to a bounded operator on $L^2$ if  
 $\Omega\in L\log L(S^{n-1})$, where 
  $L\log L(S^{n-1})$ denotes the Zygmund class of all those functions 
  $\Omega$ on $S^{n-1}$ 
  which satisfy 
  $$\int_{S^{n-1}}|\Omega(\theta)|\log(2+|\Omega(\theta)|)\, d\sigma(\theta)
  <\infty.$$  
  Furthermore, if $\Omega\in L\log L(S^{n-1})$, 
   by the method of rotations of Calder\'{o}n-Zygmund (see  \cite{2}) 
   it can be shown  that $S$ extends to a bounded operator on 
  $L^p$ for all $p\in (1,\infty)$.  
\par 
 When $n=d$ and $P(y)=y$, R. Fefferman \cite{10} 
proved  that if $h$ is bounded and $\Omega$ 
satisfies a Lipschitz condition of positive order on $S^{n-1}$, then 
the singular integral operator $T$ in (1.2)  is 
bounded on $L^p$ for $1<p<\infty$.  Namazi \cite{13} improved this result by 
replacing the Lipschitz condition  by the  condition that 
$\Omega\in L^q(S^{n-1})$ for some $q>1$.  
In \cite{7}, Duoandikoetxea and Rubio de Francia  
developed  methods which can be used to 
study mapping properties of several kinds of 
 operators in harmonic analysis including the singular integrals considered 
 in \cite{13}.  
Also, see \cite{6, 22}  for  weighted $L^p$ boundedness of singular integrals, 
and \cite{18, 19} for background materials.  
\par    
 For the rest of this note we assume that the 
  polynomial mapping $P$ in (1.2) satisfies $P(-y)=-P(y)$ and $P\neq 0$.  
We shall prove the following: 
\proclaim{Theorem 1} Let $\Omega\in L^q(S^{n-1})$, $q \in (1,2]$ and 
$h\in \Delta_s$, $s\in (1,2]$. Suppose $\Omega$ satisfies $(1.1)$. Let $T$ be as in $(1.2)$. 
Then we have 
 $$\|T(f)\|_{L^p(\Bbb R^d)}\leq C_p(q-1)^{-1}(s-1)^{-1}
 \|\Omega\|_{L^q(S^{n-1})}\|h\|_{\Delta_s}\|f\|_{L^p(\Bbb R^d)}$$ 
  for all $p\in (1,\infty)$, where the constant $C_p$ is independent of $q, s, 
  \Omega$ and $h$.  Also, the constant $C_p$  is  independent of polynomials 
  $P_j$ if we fix $\deg(P_j)$ $(j=1,2,\dots,d)$.  
\endproclaim  
In Al-Salman-Pan \cite{1}, the $L^p$ boundedness of $T$ was proved 
under the condition that $\Omega$ is a function in $L\log L(S^{n-1})$ satisfying (1.1) and $h\in \Delta_s$ for some $s>1$ (\cite{1, Theorem 1.3}).   
Also it is noted there that estimates like those in Theorem 1 
(with $s$ being fixed) 
can be used to prove the same result by applying 
an extrapolation method, 
but such estimates are yet to be proved (see \cite{1, p. 156}).  
In \cite{1}, the authors also considered singular integrals 
 defined by certain polynomial mappings $P$ which do not satisfy the condition 
 $P(-y)=-P(y)$.   
\par 
As a consequence of Theorem 1 we can give a 
different proof of \cite{1, Theorem 1.3} via an extrapolation method; in fact, we can get an improved result.   
For a positive number $a$ and a function $h$ on $\Bbb R_+$,  let
$$ 
L_a(h) = \sup_{j \in \Bbb Z} 
\int_{2^j}^{2^{j+1}}|h(r)|\left(\log(2+|h(r)|)\right)^a\, dr/r.  
$$  
We define a class $\Cal L_a$ to be the space of all those measurable functions 
$h$ on $\Bbb R_+$ which satisfy $L_a(h)<\infty$.  
 Also,  let
$$
 N_a(h)  = \sum_{m \geq 1}m^a 2^md_m(h), $$
where 
 $d_m(h) = \sup_{k \in \Bbb Z} 2^{-k}|E(k,m)|$  with 
$E(k, m) = \{r \in (2^{k}, 2^{k+1}] : 2^{m-1} < |h(r)| \leq  2^{m}\}$ 
for $m\geq 2$, 
$E(k, 1) = \{r \in (2^{k}, 2^{k+1}] : |h(r)| \leq  2 \}$.
We denote by $\Cal N_a$  the class of all  those measurable functions 
$h$ on $\Bbb R_+$ such that  $N_a(h)<\infty$.  
Then we readily see that $N_a(h) < \infty$ implies 
$L_a(h)< \infty$.   
Conversely, 
if $L_{a+b}(h) < \infty$ for some $b > 1$, then $N_a(h) < \infty$.  
To see this, note that 
$$2^m m^{a+b} 2^{-k}|E(k,m)|  \leq C\int_{E(k,m)}
|h(r)|\left(\log(2+|h(r)|)\right)^{a+b}\, dr/r \leq CL_{a+b}(h)$$  
for $m\geq 2$; thus 
$ N_a(h) \leq 2d_1(h) + CL_{a+b}(h)\sum_{m \geq 2} m^{-b} < \infty$.    
By Theorem 1 and an extrapolation method we have the following: 
\proclaim{Theorem 2}  Suppose $\Omega$ is a function in $L\log L(S^{n-1})$ 
satisfying $(1.1)$ and $h\in \Cal N_1$.  Let $T$ be as in $(1.2)$. 
Then
 $$\|T(f)\|_{L^p(\Bbb R^d)}\leq C_p\|f\|_{L^p(\Bbb R^d)}$$ 
  for all $p\in (1,\infty)$, where $C_p$ is independent of polynomials $P_j$ 
  if the polynomials are of fixed degree.   
\endproclaim
By Theorem 2  and the remark preceding it we see that 
$T$ is bounded on $L^p$ for all $p\in (1,\infty)$ if  
$\Omega$ is as in Theorem 2 and $h\in \Cal L_{a}$ for some $a>2$.  
\par 
When $n=d$, $P(y)=y$, $\Omega$ is as in Theorem 2   
and $h$ is a constant function, it is known that $T$ is of weak type $(1,1)$; 
see \cite{5, 15}.     
Also,  see  \cite{4, 9, 11, 12, 16, 20, 21} for related results.  
\par 
 In Section 2, we shall prove Theorem 1.   
  Applying the methods of \cite{7} 
 involving the Littlewood-Paley theory and  using results of 
  \cite{8, 14},   
  we shall prove  $L^p$ estimates for  certain maximal and 
  singular integral operators related to the operator $T$ in Theorem 1
   (Lemmas 1 and 2).  Lemma 1 is used to prove Lemma 2.  By Lemma 2 we can 
   easily prove Theorem 1.  A key idea of 
  the proof of Theorem 1 is  to apply a Littlewood-Paley 
  decomposition adapted to a 
  suitable  lacunary sequence depending on $q$ and $s$ for which 
  $\Omega\in L^q(S^{n-1})$ and $h\in \Delta_s$.  
  The method of  appropriately choosing the  lacunary sequence 
  was inspired by \cite{1}, where, 
  in a somewhat different way from ours, a similar method was used to study  
    several classes of singular integrals.    
  \par 
  We shall prove Theorem 2 in Section 3.  
  Finally, in Section 4, we consider the maximal operator   
  $$T^*(f)(x) 
= \sup_{N, \epsilon > 0}\left|\int_{\epsilon<|y|<N}
f(x - P(y))K(y)\, dy\right|,   \tag 1.3 $$ 
 where $P$ and $K$ are as in (1.2).  
We shall prove  analogs of Theorems 1 and 2 for the operator $T^*$. 
 \par 
    Throughout this note,  
 the letter $C$ will be used to denote non-negative constants which may be 
different in different occurrences.

\head 2. Proof of Theorem 1  \endhead  

Let $\Omega, h$ be as in Theorem 1.  
We consider the singular integral $T(f)$ defined in (1.2).   
Let $\rho\geq 2$ and    
$E_k=\{x\in \Bbb R^n: \rho^k<|x| \leq \rho^{k+1}\}$. 
  Then    
$T(f)(x) = \sum_{-\infty}^{\infty}\sigma_k * f(x)$,      
where $\{\sigma_k\}$ is a sequence of Borel measures on $\Bbb R^d$ such that 
$$\sigma_k*f(x)=\int_{E_k}f(x-P(y))K(y)\,dy. \tag 2.1 $$ 
We note that 
$$(\sigma_k*f)\hat{\phantom{t}}(\xi)=\hat{f}(\xi)\int_{E_k}  
e^{-2\pi i\langle P(y), \xi\rangle}K(y)\,dy. $$  
\par 
We write 
$$P(y)=\sum_{j=1}^\ell Q_j(y), \qquad  Q_j(y)=\sum_{|\gamma|=N(j)}a_\gamma 
y^\gamma \quad (a_\gamma\in \Bbb R^d),$$ 
where $Q_j\neq 0$, 
$1\leq N(1)<N(2)<\dots<N(\ell)$, $\gamma=(\gamma_1,\dots,\gamma_n)$ is a 
multi-index, $y^\gamma=y_1^{\gamma_1}\dots y_n^{\gamma_n}$ and 
 $|\gamma|=\gamma_1+\dots + \gamma_n$.   
 Let  $\beta_m=\rho^{N(m)}$ and 
$\alpha_m=(q - 1)(s - 1)/(2qsN(m))$ for $1\leq m \leq \ell $. 
 Put $P^{(m)}(y)=\sum_{j=1}^mQ_j(y)$ and  
define a sequence $\mu^{(m)}=\{\mu_k^{(m)}\}$ of positive measures on 
$\Bbb R^d$  by 
$$\mu_k^{(m)}*f(x)=\int_{E_k} f\left(x-P^{(m)}(y)\right)|K(y)| \,dy$$  
for $m=1,2,\dots, \ell$.  Also, define 
$\mu^{(0)}=\{\mu_k^{(0)}\}$ by  $\mu_k^{(0)}= 
(\int_{E_k}|K(y)| \,dy)\delta$, where $\delta$ is Dirac's delta function 
on $\Bbb R^d$.  
For a sequence $\nu=\{\nu_k\}$ of finite Borel measures on $\Bbb R^d$, 
we define the  
maximal operator $\nu^*$ by $\nu^*(f)(x) = \sup_k ||\nu_k|*f(x)|$, 
where   $|\nu_k |$ denotes the total variation.   
We consider the maximal operators $\left(\mu^{(m)}\right)^*$   
($0\leq m\leq \ell$).   
We also write $\left(\mu^{(\ell)}\right)^*=\mu^*_\rho$.       
\par
 Let 
$$L_j(\xi)= (\langle a_{\gamma(j,1)},\xi\rangle, 
\langle a_{\gamma(j,2)},\xi\rangle, \dots, \langle a_{\gamma(j,r_j)},
\xi\rangle), 
$$  
where $\{\gamma(j,k)\}_{k=1}^{r_j}$ is an enumeration of 
$\{\gamma\}_{|\gamma|=N(j)}$ for  $1\leq j \leq \ell $.  
Then $L_j$ is a linear mapping from $\Bbb R^d$ to $\Bbb R^{r_j}$. 
Let $s_j=\text{\rm rank}\, L_j$. There exist non-singular 
linear transformations $R_j:\Bbb R^d\to \Bbb R^d$ and $H_j:\Bbb R^{s_j}\to 
\Bbb R^{s_j}$ such that  
$$|H_j\pi^d_{s_j}R_j(\xi)|\leq |L_j(\xi)|\leq C|H_j\pi^d_{s_j}R_j(\xi)|,  $$ 
 where $\pi^d_{s_j}(\xi)=(\xi_1,\dots, \xi_{s_j})$ is the projection and 
  $C$ depends only on $r_j$ (a proof can be found in \cite{8}).   
Let $\{\sigma_k^{(m)}\}$ ($0\leq m\leq\ell$) be a sequence of
 Borel measures on $\Bbb R^d$ such that  
$$\sigma_k^{(m)}*f(x)=\int_{E_k} f\left(x-P^{(m)}(y)\right)K(y)\,dy $$ 
for  $m=1, 2, \dots, \ell$, while  $\sigma_k^{(0)}=0$. 
Let $\varphi\in C_0^{\infty}(\Bbb R)$ be supported in $\{|r| \leq 1 \}$ 
  and   $\varphi(r) = 1$ for $|r|<1/2$.    
Define a sequence $\tau^{(m)}=\{\tau_k^{(m)}\}$ of Borel measures by
$$\hat{\tau}_k^{(m)}(\xi)= 
\hat{\sigma}_k^{(m)}(\xi)\Phi_{k,m}(\xi)
-\hat{\sigma}_k^{(m-1)}(\xi)\Phi_{k,m-1}(\xi)   \tag 2.2 $$ 
 for $m=1, 2, \dots, \ell$, where 
 $$\Phi_{k,m}(\xi)=\prod_{j=m+1}^\ell
  \varphi\left(\beta_j^{k}|H_j\pi^d_{s_j}R_j(\xi)|\right)$$
if $0\leq m\leq \ell-1$ and $\Phi_{k,\ell}=1$.  
Then $\sigma_k=\sigma_k^{(\ell)}=\sum_{m=1}^\ell \tau_k^{(m)}$.    
We note that 
$$\Phi_{k,m}(\xi)\varphi\left(\beta_m^{k}|H_m\pi^d_{s_m}R_m(\xi)|\right)
=\Phi_{k,m-1}(\xi) \quad (1\leq m\leq \ell).  $$ 
For $1\leq m\leq \ell$, let  
$T_\rho^{(m)}(f)=\sum_k \tau_k^{(m)}*f$.  
Then $T=\sum_{m=1}^\ell T_\rho^{(m)}$.  
\par 
For $p\in (1,\infty)$ we put  $p'=p/(p-1)$ and 
 $\delta (p) = |1/p - 1/p'|$.   
Let $\theta  \in (0, 1)$.     
Then we have the following $L^p$  
 estimates for  $(\mu^{(m)})^*$ and $T_\rho^{(m)}$.      
\proclaim{Lemma 1} For  $p >1+\theta$ and $0\leq j\leq \ell$, we have 
$$\left\|(\mu^{(j)})^*(f)\right\|_{L^p(\Bbb R^d)} \leq C 
(\log \rho)\|\Omega\|_{L^q(S^{n-1})}\|h\|_{\Delta_s}
\left(1-\rho^{-\theta/(2q's')}\right)^{-2/p}
\|f\|_{L^p(\Bbb R^d)}.      \tag 2.3  $$  
\endproclaim 
\proclaim{Lemma 2} For $p \in (1+\theta, (1+\theta)/\theta)$ and     
 $1\leq m\leq \ell$,  
we have 
$$\|T_\rho^{(m)}(f)\|_{L^p(\Bbb R^d)} 
\leq C (\log \rho)\|\Omega\|_{L^q(S^{n-1})}\|h\|_{\Delta_s}
\left(1-\rho^{-\theta/(2q's')}\right)^{-1 -\delta (p)}
\|f\|_{L^p(\Bbb R^d)}.  $$ 
\endproclaim  
The constants  $C$ in Lemmas 1 and 2 are  
independent of  
$q, s\in (1,2]$,  $\Omega\in L^q(S^{n-1})$, 
$h\in \Delta_s$,   $\rho$  
and the coefficients of the polynomials  $P_k$ $(1\leq k\leq d)$. 
\par 
We prove Lemma 2 first, taking Lemma 1 for granted for the moment.  Let 
$A=(\log \rho)\|\Omega\|_{L^q(S^{n-1})}\|h\|_{\Delta_s}$ 
and $B=\left(1-\beta_m^{-\theta\alpha_m}\right)^{-1}=
\left(1-\rho^{-\theta/(2q's')}\right)^{-1}$.  
Then, we have the following estimates:  
$$ \|\tau_k^{(m)}\| \leq c_1A \quad 
(\|\tau_k^{(m)}\|=|\tau_k^{(m)}|(\Bbb R^d)),      
\tag 2.4 $$
$$|\hat{\tau}_k^{(m)}(\xi)| \leq c_2A
\left(\beta_m^{k}|L_m(\xi)|\right)^{-\alpha_m},   \tag 2.5 $$ 
$$|\hat{\tau}_k^{(m)}(\xi)|\leq c_3A\left(\beta_m^{k+1}|L_m(\xi)|
\right)^{\alpha_m},   \tag 2.6 $$  
$$
\left\|(\tau^{(m)})^*(f)\right\|_p\leq C_pAB^{2/p}\|f\|_p  
\quad \text{for $p>1+\theta$},  
\tag 2.7   
$$ 
for some  constants $c_i$ ($1\leq i\leq 3$) and  $C_p$,  
where we simply write $\|f\|_{L^p(\Bbb R^d)}=\|f\|_p$.  
\par 
Now we prove the estimates (2.4)--(2.7). 
First we see that
$$\align 
 \|\tau_k^{(m)}\| &\leq C\left(\|\sigma^{(m)}_k\|+\|\sigma^{(m-1)}_k\|\right) 
 \tag 2.8  
 \\
 &\leq C\|\Omega\|_1\int_{\rho^k}^{\rho^{k+1}}|h(r)|\, dr/r 
 \leq C(\log \rho)\|\Omega\|_1\|h\|_{\Delta_1}. 
 \endalign $$ 
From this (2.4) follows. 
 To prove (2.5),  define 
$$F(r,\xi)=\int_{S^{n-1}}\Omega(\theta)
\exp(-2\pi i\langle\xi, P^{(m)}(r\theta)\rangle)\, d\sigma(\theta).     
$$ 
Then,  via H\"older's inequality,  for $s\in (1,2]$ we see that  
$$\align  
&|\hat{\sigma}_k^{(m)}(\xi)| 
= \left| \int_{\rho^{k}}^{\rho^{k+1}}h(r) F(r,\xi)\, dr/r \right|    
 \tag 2.9 
\\
&\leq \left(\int_{\rho^{k}}^{\rho^{k+1}}|h(r)|^s\,dr/r \right)^{1/s}
\left(\int_{\rho^{k}}^{\rho^{k+1}}\left|F(r,\xi)\right|^{s'}\,
dr/r\right)^{1/s'}\\ 
&\leq C(\log \rho)^{1/s}\|h\|_{\Delta_s}
\|\Omega\|_1^{(s'-2)/s'}\left(\int_{\rho^{k}}^{\rho^{k+1}}
\left|F(r,\xi)\right|^2\, dr/r\right)^{1/s'}. 
\endalign$$  
We need the following estimates for the last integral:  

\proclaim{Lemma 3}  Let $1<q\leq 2$ and $\Omega\in L^q(S^{n-1})$.  
Then there exists  a constant $C>0$ independent of $q, \rho, \Omega$ and 
 the coefficients of the polynomial components of  $P^{(m)}$  such that 
$$\int_{\rho^{k}}^{\rho^{k+1}}\left|F(r,\xi)\right|^2\, dr/r 
\leq C(\log \rho)\left(\beta_m^k|L_m(\xi)|\right)^{-1/(2q'N(m))}\|\Omega\|_q^2. $$
\endproclaim  
\demo{Proof}
Take an integer $\nu$ such that 
$2^\nu<\rho\leq 2^{\nu+1}$.  By the proof of Proposition 5.1 of \cite{8}
 we have 
$$\align 
\int_{\rho^{k}}^{\rho^{k+1}}\left|F(r,\xi)\right|^2\, dr/r  
&=\int_{1}^{\rho}\left|F(\rho^k r,\xi)\right|^2\, dr/r  \leq \sum_{j=0}^\nu   
\int_{2^j}^{2^{j+1}}\left|F(\rho^k r,\xi)\right|^2\, dr/r  
\\  
&\leq \sum_{j=0}^\nu   
C\left(\int_{1}^{2}\left|F(2^j\rho^k r,\xi)\right|^{q'}\, dr/r\right)^{2/q'}
\\   
&\leq  \sum_{j=0}^\nu C\left(2^{jN(m)}\rho^{kN(m)}|L_m(\xi)|
\right)^{-1/(2N(m)q')}\|\Omega\|_q^2 
\\
&\leq C(\log\rho)\left(\rho^{kN(m)}|L_m(\xi)|
\right)^{-1/(2N(m)q')}\|\Omega\|_q^2.     
\endalign $$ 
This completes the proof of Lemma 3.     
\enddemo 
By  (2.9) and  Lemma 3 we have $|\hat{\sigma}_k^{(m)}(\xi)| \leq CA
\left(\beta_m^{k}|L_m(\xi)|\right)^{-\alpha_m}$.  
Also, we have $\|\sigma^{(m-1)}_k\| \leq CA$ by (2.8).  
 We can prove the estimate (2.5) by using 
 these estimates in the definition of $\tau_k^{(m)}$ in (2.2) and by  
 noting that $\varphi$ is compactly supported.    
Next, to prove (2.6), using (1.1) when $m=1$,  we see that  
 $$\align 
 &|\hat{\tau}_k^{(m)}(\xi)|
\leq \left|\left(\hat{\sigma}_k^{(m)}(\xi)-\hat{\sigma}_k^{(m-1)}(\xi)
\right)\Phi_{k,m}(\xi)\right| 
+\left|\left(\Phi_{k,m}(\xi)-\Phi_{k,m-1}(\xi)\right)
\hat{\sigma}_k^{(m-1)}(\xi)\right| 
 \\ 
 &\leq C\|\Omega\|_1 \beta_m^{k+1}|L_m(\xi)|
\int_{\rho^k}^{\rho^{k+1}}|h(r)|\, dr/r  
+C\|\sigma_k^{(m-1)}\| \beta_m^{k}|L_m(\xi)|
\\ 
   &\leq C(\log\rho)
\|\Omega\|_1\|h\|_{\Delta_1}\beta_m^{k+1}|L_m(\xi)|,  
\endalign $$ 
where to get the last inequality we have used (2.8).       
By this and (2.8),   we have  
$$|\hat{\tau}_k^{(m)}(\xi)|\leq C(\log\rho) 
\|\Omega\|_1 \|h\|_{\Delta_1}\left(\beta_m^{k+1}|L_m(\xi)|
\right)^{c}  $$  
for all $c\in (0,1]$, which implies (2.6). 
Finally, the estimate (2.7) follows from 
 Lemma 1  since 
$$ 
\left\|(\tau^{(m)})^*(f)\right\|_p 
\leq C\left\|(\mu^{(m)})^*(|f|)\right\|_p+
C\left\|(\mu^{(m-1)})^*(|f|)\right\|_p                              
\leq CAB^{2/p}\|f\|_p  $$  
for $p>1+\theta$,  where  the first inequality can be seen by change of 
variables and a well-known result on maximal functions (see \cite{8}).  
\par 
Let $\{\psi_k\}_{- \infty}^{\infty}$ be a sequence of functions  in
 $C^{\infty}((0, \infty))$  such that 
$$ \text{supp} (\psi_k) \subset [\beta_m^{-k-1}, \beta_m^{-k+1}], \quad
\sum_k \psi_k(t)^2 = 1, \quad |(d/dt)^j\psi_k (t)| \leq c_j/t^j \ \  
(j=1,2, \dots),  $$ 
where the constants $c_j$ are independent of $\beta_m$.  
 Define an operator $S_k$ by
$ \left(S_k(f)\right)\hat{\phantom{t}}(\xi) = 
\psi_k\left(|H_m\pi^d_{s_m}R_m(\xi)|\right)\hat{f}(\xi)$ and let 
$$V_j^{(m)}(f) =  \sum_{-\infty}^{\infty} S_{j+k}\left(\tau_k^{(m)} 
* S_{j+k}(f)\right).$$
Then by Plancherel's theorem and the estimates  (2.4)--(2.6)  we have
$$\align
\left\|V_j^{(m)}(f) \right\|_2^2 
&\leq \sum_k C\int_{D(j + k)} |\hat{\tau}_k^{(m)}(\xi)|^2|\hat{f}(\xi)|^2 
\, d\xi \tag 2.10 \\
&\leq CA^2\min\left(1, \beta_m^{-2\alpha_m(|j|-2)}\right) \sum_k  
\int_{D(j + k)} |\hat{f}(\xi)|^2 \, d\xi \\
&\leq CA^2\min\left(1, \beta_m^{-2\alpha_m(|j|-2)}\right)\|f\|^2_2,
\endalign$$
where $D(k) = \{\beta_m^{-k-1} \leq |H_m\pi^d_{s_m}R_m(\xi)| 
\leq \beta_m^{-k+1}\}$.  
\par 
Applying the proof of Lemma in \cite {7, p. 544} and using the estimates 
(2.4) and (2.7), we 
can prove the following.  
\proclaim{Lemma 4}   
Let  $u\in (1+\theta, 2]$. 
Define a number $v$ by $1/v-1/2=1/(2u)$. 
Then we have the vector valued inequality   
$$
\left\|\left(\sum |\tau_k^{(m)}*g_k|^2\right)^{1/2}\right\|_{v} 
\leq (c_1C_u)^{1/2}AB^{1/u}
\left\|\left(\sum |g_k|^2\right)^{1/2}\right\|_{v},   \tag 2.11 
 $$  
 where the constants $c_1$ and $C_u$ are as in $(2.4)$ and $(2.7)$, 
 respectively.     
\endproclaim 
\par 
By the Littlewood-Paley theory  we have 
$$\gather 
\|V_j^{(m)}(f)\|_{p} 
\leq c_{p}\left\| \left(\sum_k |\tau_k^{(m)}*S_{j+k}(f)|^2\right)^{1/2}
\right\|_{p},  \tag 2.12 
\\ 
\left\|\left(\sum_k |S_{k}(f)|^2\right)^{1/2}\right\|_{p}\leq c_p\|f\|_{p}, 
\tag 2.13  
\endgather 
$$ 
where $1<p<\infty$ and $c_p$ is independent of $\beta_m$ and the linear 
transformations $R_m, H_m$.  
Suppose that $1+\theta<p\leq 4/(3-\theta)$.  
Then we can find $u\in (1+\theta, 2]$ such that $1/p=1/2+(1-\theta)/(2u)$. 
Let $v$ be defined by $u$ as in Lemma 4. Then 
by  (2.11)--(2.13) we have 
$$
\|V_j^{(m)}(f)\|_{v} \leq CAB^{1/u}\|f\|_{v}.   \tag 2.14   
$$  
 Since $1/p=\theta/2+(1-\theta)/v$, interpolating between (2.10) and (2.14), 
 we have 
 $$ \|V_j^{(m)}(f)\|_p
 \leq CAB^{(1-\theta)/u}\min\left(1, \beta_m^{-\theta\alpha_m(|j|-2)}\right)
 \|f\|_p.  $$  
 It follows that 
$$\align 
\|T_\rho^{(m)}(f)\|_{p} &\leq \sum_j \|V_j^{(m)}(f)\|_{p} 
\leq CAB^{(1 - \theta)/u}
(1-\beta_m^{-\theta \alpha_m})^{-1}\|f\|_{p}  \tag 2.15 
\\ 
&\leq C AB^{2/p}\|f\|_{p}, 
\endalign $$    
where we have used the inequality 
$\sum \min\left(1, \beta_m^{-\theta\alpha_m(|j|-2)}\right)\leq 
5\left(1 - \beta_m^{- \theta\alpha_m}\right)^{-1}$.   
We also have $\|T_\rho^{(m)}(f)\|_{2}\leq \sum \|V_j^{(m)}(f)\|_{2}
\leq C AB\|f\|_{2}$ by (2.10), 
 since $B\geq \left(1-\beta_m^{-\alpha_m}\right)^{-1}$. 
By duality and interpolation,  we can now get the conclusion of Lemma 2.  
\par  
Next, we give a proof of  Lemma 1. We prove Lemma 1 by induction on $j$.   
 Now we  assume (2.3) for $j=m-1$, $1\leq m\leq \ell$,  and prove (2.3) for 
 $j=m$.    Let $\varphi\in C_0^{\infty}(\Bbb R)$ be as above.  
  Define a sequence $\eta^{(m)}=\{\eta_k^{(m)}\}$ of Borel measures  
  on $\Bbb R^d$ by 
$$\hat{\eta}_k^{(m)}(\xi)=  
\varphi\left(\beta_m^{k}|H_m\pi^d_{s_m}R_m(\xi)|\right)\hat{\mu}_k^{(m-1)}(\xi).$$                   
  Then by (2.3) with $j=m-1$, we have 
  $$
  \left\|(\eta^{(m)})^*(f)\right\|_p 
  \leq C\left\|(\mu^{(m-1)})^*(f)\right\|_p  
  \leq CAB^{2/p}\|f\|_p    \tag 2.16
$$ 
for $p>1+\theta$. 
Furthermore, we have the following: 
$$\align 
  \|\eta_k^{(m)}\|+\|\mu_k^{(m)}\| &\leq C\|\mu_k^{(m-1)}\|+\|\mu_k^{(m)}\| 
\leq C\|\Omega\|_1\int_{\rho^k}^{\rho^{k+1}}|h(r)|\, dr/r   \tag 2.17 
\\ 
&\leq C(\log \rho)\|\Omega\|_1\|h\|_{\Delta_1}\leq CA,     
\endalign $$ 
$$\align 
|\hat{\mu}_k^{(m)}(\xi)-\hat{\eta}_k^{(m)}(\xi)|&\leq C(\log\rho) 
\|\Omega\|_1 \|h\|_{\Delta_1}
\left(\beta_m^{k+1}|L_m(\xi)|\right)^{\alpha_m}   \tag 2.18
\\ 
&\leq CA\left(\beta_m^{k+1}|L_m(\xi)|\right)^{\alpha_m},     
\endalign$$
 $$|\hat{\mu}_k^{(m)}(\xi)| \leq CA
\left(\beta_m^{k}|L_m(\xi)|\right)^{-\alpha_m},   \tag 2.19
$$ 
$$|\hat{\eta}_k^{(m)}(\xi)| \leq C(\log \rho)\|h\|_{\Delta_1}\|\Omega\|_1
\left(\beta_m^{k}|L_m(\xi)|\right)^{-\alpha_m}
\leq CA\left(\beta_m^{k}|L_m(\xi)|\right)^{-\alpha_m}.  \tag 2.20 $$  
\par 
To see  (2.18) we note that    
$$ 
 |\hat{\mu}_k^{(m)}(\xi)-\hat{\eta}_k^{(m)}(\xi)|  \leq 
 |\hat{\mu}_k^{(m)}(\xi)-\hat{\mu}_k^{(m-1)}(\xi)|+\left|\left(\varphi
 \left(\beta_m^{k}|H_m\pi^d_{s_m}R_m(\xi)|\right)-1\right)
 \hat{\mu}_k^{(m-1)}(\xi)\right|. 
$$ 
Thus arguing as in the proof of (2.6), 
we have the first inequality of (2.18). 
 The estimate (2.19) follows from the arguments used to prove (2.5). 
Also, we can see the first inequality of (2.20) by the definition  
of $\eta_k^{(m)}$ and (2.17). 
 \par  
Since $\|(\mu^{(m)})^*(f)\|_{\infty}\leq CA\|f\|_\infty$, 
by taking into account an interpolation, 
it suffices to prove (2.3) with $j=m$ for  $p\in (1+\theta, 2]$.  
 Define a sequence $\nu^{(m)}=\{\nu_k^{(m)}\}$ of
Borel measures by $\nu_k^{(m)} = \mu_k^{(m)} - \eta_k^{(m)}$.   
Let
$$g_m(f)(x) = \left( \sum \left|\nu_k^{(m)} *f(x) \right |^2 \right)^{1/2}.$$
Then 
$$(\mu^{(m)})^*(f) \leq g_m(f) + (\eta^{(m)})^*(|f|). 
 \tag 2.21 $$
 Thus, by (2.16), to get  (2.3) with $j=m$ 
 it suffices to prove $\|g_m(f)\|_p \leq CAB^{2/p}\|f\|_p $ 
for $p\in (1+\theta, 2]$ with an appropriate constant $C$.      
By a well-known property of Rademacher's functions, this follows from  
$$ \left \|U^{(m)}_\epsilon(f) \right \|_p \leq C AB^{2/p}
\|f\|_p \tag 2.22 $$
for  $p\in (1+\theta, 2]$, where 
$ U^{(m)}_\epsilon(f) = \sum_k \epsilon_k \nu^{(m)}_k * f$  
with $\epsilon = \{\epsilon_k \}$, $\epsilon_k = 1$ or $-1$,   
 and the constant $C$ is independent of $\epsilon$. 
\par
The estimate (2.22) is a consequence of the following: 
\proclaim{Lemma 5} 
We define a sequence $\{p_j \}_1^{\infty}$ by $p_1 = 2$ and  
$1/p_{j+1} = 1/2  + (1 - \theta)/(2p_j)$ for $j \geq 1$.  
$($We note that $1/p_j=(1-a^j)/(1+\theta)$,  where $a=(1-\theta)/2$, so  
$\{p_j\}$ is decreasing and converges to $1 + \theta$.$)$ 
  Then, for $j \geq 1$ we have  
$$\left \|U^{(m)}_{\epsilon}(f)\right \|_{p_j} 
\leq C_j AB^{2/p_j}\left \|f\right \|_{p_j}. $$
\endproclaim 
\demo{Proof} 
Let 
$$U^{(m)}_j(f) =  \sum_{k=-\infty}^{\infty}\epsilon_k S_{j+k}\left(\nu_k^{(m)} 
* S_{j+k}(f)\right).$$
Then by Plancherel's theorem and the estimates (2.17)--(2.20), as in (2.10) 
  we have
$$
\left\|U^{(m)}_j(f) \right\|_2 
\leq CA\min\left(1, \beta_m^{-\alpha_m(|j|-2)}\right)\|f\|_2.   \tag 2.23  
$$  
It follows that 
$
\left\|U^{(m)}_{\epsilon}(f) \right\|_2 \leq \sum_j\|U^{(m)}_j(f)\|_2 
\leq CAB\|f\|_2$.   
If we denote by $A(s)$ the assertion of Lemma 5 for $j=s$, this proves  
 $A(1)$.  
\par  
Now we derive $A(s+1)$ from $A(s)$ assuming that $A(s)$ holds,  
which will complete the proof of  Lemma 5 by induction.  
Using (2.21), we see that 
$$ (\nu^{(m)})^*(f) \leq (\mu^{(m)})^*(|f|) + 
(\eta^{(m)})^*(|f|) \leq g_m(|f|) + 2(\eta^{(m)})^*(|f|). 
$$
 Note that $A(s)$ implies  
$\|g_m(f)\|_{p_s} \leq CAB^{2/p_s}\|f\|_{p_s}$.  
 By this and  (2.16) we have
$$\left\|(\nu^{(m)})^*(f)\right\|_{p_s} \leq \|g_m(|f|)\|_{p_s} 
+ 2\left\|(\eta^{(m)})^*(|f|)\right\|_{p_s} 
\leq CAB^{2/p_s}\|f\|_{p_s}. \tag 2.24 $$
By (2.17),  (2.23) and (2.24),  we can now apply the arguments used in 
the proof of (2.15) to get 
 $A(s + 1)$. This  completes the proof of Lemma 5.
\enddemo
Now we prove (2.22) for $p\in (1+\theta,2]$.   Let $\{p_j\}_1^{\infty}$ be 
as in Lemma 5.  Then we have $p_{N+1} < p \leq p_N$ 
for some $N$. 
Thus, interpolating between the estimates of Lemma 5 for $j=N$ and $j=N+1$, 
  we have (2.22).  
This proves (2.3) for  $j=m$. 
 \par 
Finally, we can easily see that
$(\mu^{(0)})^*(f) \leq C(\log\rho)\|\Omega\|_1\|h\|_{\Delta_1}|f|$ 
(see (2.17)), which implies the estimate  (2.3) for $j=0$.    
Therefore,  by induction we have (2.3) for all $0\leq j\leq \ell$.   
This completes the proof of Lemma 1. 
\par 
Now we can prove Theorem 1. 
Since $\theta\in (0,1)$ is  arbitrary,
by taking $\rho=2^{q's'}$ in Lemma 2  we have
$$\|T_{2^{q's'}}^{(m)}(f)\|_{p}\leq C_p(q-1)^{-1}(s-1)^{-1}
\|\Omega\|_{q}\|h\|_{\Delta_s}\|f\|_{p}$$ 
 for all $p\in (1,\infty)$. This 
completes the proof of Theorem 1, since $T=\sum_{m=1}^\ell T_\rho^{(m)}$.  

\head 3. Proof of Theorem 2  \endhead  

 Theorem 2 can be proved  by  Theorem 1 and an extrapolation 
 argument. 
 Let $T(f)$ be the singular integral in (1.2). 
 We also write $T(f)=T_{h,\Omega}(f)$. 
We fix $q\in (1,2]$, $\Omega\in L^q(S^{n-1})$,  $p \in (1, \infty)$ and 
a function $f$ with $\|f\|_p\leq 1$ and put     
 $S(h) = \|T_{h,\Omega}(f)\|_p$.  
Then we have the following subadditivity:  
$$S(h + k ) \leq S(h) + S(k). \tag 3.1 $$
Set $E_1 = \{r \in \Bbb R_+ : |h(r)| \leq 2 \}$ and 
$E_m = \{r \in \Bbb R_+ : 2^{m-1} < |h(r)| \leq 2^{m}\}$  
for $m \geq 2$. Then, applying Theorem 1, we see that  
$$S\left(h\chi_{E_m}\right) 
\leq C(q - 1)^{-1}(s - 1)^{-1}
\|\Omega\|_q\|h\chi_{E_m}\|_{\Delta_s} \tag 3.2$$
for $s \in (1, 2]$, where $\chi_{E}$ 
denotes  the characteristic function of a  set $E$.    
Now we follow the extrapolation argument of Zygmund \cite 
{23, Chap. XII, pp. 119--120}.  
First, note that
$$\|h\chi_{E_m}\|_{\Delta_{1+1/m}}  \leq 2^{m}d_m^{m/(m+1)}(h) $$
for $ m \geq 1$, where $d_m(h)$ is as in Section 1.  Using this and  
 (3.2) we see that 
$$
\align
\sum_{m \geq 1} S\left(h\chi_{E_m}\right) 
&\leq C(q - 1)^{-1}\|\Omega\|_q \sum_{m \geq 1}
m\|h\chi_{E_m}\|_{\Delta_{1+1/m}} 
\\
&\leq C(q - 1)^{-1}\|\Omega\|_q \sum_{m \geq 1}
m2^{m}d_m^{m/(m+1)}(h). 
\endalign
$$
Recalling the definition of $N_a(h)$, we have
$$\align 
&\sum_{m \geq 1}m2^{m}d_m^{m/(m+1)}(h) = 
\sum_{d_m(h) < 3^{-m}}m2^{m}d_m^{m/(m+1)}(h)  
+  \sum_{d_m(h) \geq 3^{-m}}m2^{m}d_m^{m/(m+1)}(h) \\
&\leq \sum_{m \geq 1}m2^{m}3^{-m^2/(m+1)} 
+ \sum_{m \geq 1}m2^{m}d_m(h)3^{m/(m+1)} \leq C(1 + N_1(h)).
\endalign$$
Therefore, by (3.1) we see that 
$$S(h)\leq \sum_{m \geq 1} S\left(h\chi_{E_m}\right) \leq 
C (q - 1)^{-1}\|\Omega\|_q\left(1 + N_1(h)\right).  \tag 3.3 
$$ 
\par 
Next, fix $h\in \Cal N_1$, $p\in(1,\infty)$ and $f$ with $\|f\|_p\leq 1$ 
and let $R(\Omega)=\|T_{h,\Omega}(f)\|_p$.  Put $e_m=\sigma(F_m)$ for 
$m\geq 1$, where 
$F_m=\{\theta\in S^{n-1}: 
2^{m-1}<|\Omega(\theta)|\leq 2^m\}$ for $m\geq 2$ and 
$F_1=\{\theta\in S^{n-1}: |\Omega(\theta)|\leq 2\}$.  
We decompose $\Omega$ as $\Omega=\sum_{m=1}^\infty \Omega_m$, where 
$\Omega_m=\Omega\chi_{F_m}-\sigma(S^{n-1})^{-1}\int_{F_m}\Omega\,d\sigma$. 
We note that $\int\Omega_m\,d\sigma=0$, 
$\|\Omega_m\|_r\leq C2^m e_m^{1/r}$ for $1<r<\infty$. 
Now, by (3.3) and the subadditivity of $R(\Omega)$ we see that 
$$\align
R(\Omega)&\leq \sum_{m\geq 1}R(\Omega_m) 
\leq C \left(1 + N_1(h)\right)\sum_{m\geq 1}m\|\Omega_m\|_{1+1/m} 
\\ 
&\leq C \left(1 + N_1(h)\right)\sum_{m\geq 1}m 2^m e_m^{m/(m+1)} 
= C\left(1 + N_1(h)\right)\left(\sum_{e_m < 3^{-m}} + \sum_{e_m \geq 3^{-m}}
\right)  
\\
&\leq C \left(1 + N_1(h)\right)\left(\sum_{m \geq 1}m2^m3^{-m^2/(m+1)}
 + \sum_{m \geq 1}m2^me_m3^{m/(m+1)}\right)
 \\
&\leq C \left(1 + N_1(h))\right)
\left(1+\int_{S^{n-1}}|\Omega(\theta)|\log(2+|\Omega(\theta)|)\, 
d\sigma(\theta)\right).
\endalign$$   
This completes the proof of Theorem 2.  
   
\head  4. Estimates for maximal functions \endhead 

 For the maximal operator $T^*$ in $(1.3)$ 
 we have a result similar to Theorem 1.  
\proclaim{Theorem 3} Let $q \in (1,2]$, $s\in (1,2]$ and $\Omega\in 
L^q(S^{n-1})$, $h\in \Delta_s$.  
Suppose  $\Omega$ satisfies $(1.1)$.  
Then  we have 
 $$\|T^*(f)\|_{L^p(\Bbb R^d)}\leq C_p(q-1)^{-1}(s-1)^{-1}
 \|\Omega\|_{L^q(S^{n-1})}\|h\|_{\Delta_s}\|f\|_{L^p(\Bbb R^d)}$$ 
  for all $p\in (1,\infty)$, where $C_p$ is independent of $q$, $s$, $\Omega$ 
  and $h$. 
\endproclaim  
As Theorem 1 implies Theorem 2, 
we have the following as a consequence of Theorem 3.  
\proclaim{Theorem 4} Let $\Omega$ be a function in $L\log L(S^{n-1})$ 
satisfying $(1.1)$ and $h\in \Cal N_1$. 
  Then 
  $$\|T^*(f)\|_{L^p(\Bbb R^d)}\leq C_p\|f\|_{L^p(\Bbb R^d)}$$ 
  for all $p\in (1,\infty)$. 
  
\endproclaim  
As in the cases of Theorems 1 and 2, the constants $C_p$ of Theorems 3 and 4 
are also independent of polynomials $P_j$ if we fix $\deg(P_j)$ 
($j=1,2,\dots,d$).  
When $\Omega$ is as in Theorem 4 and  $h\in \Delta_s$ 
for some $s>1$,   the $L^p$ boundedness of $T^*$ was proved in 
\cite{1}. 
When $n=d$, $P(y)=y$, $\Omega\in L^q$ for some $q>1$ and $h$ is bounded, 
the $L^p$ boundedness of $T^*$ is due to \cite{3}. 
\par 
We use the following to prove Theorem 3. 
\proclaim{Lemma 6} 
Let $\tau^{(m)}=\{\tau^{(m)}_k\}$ $(1\leq m\leq \ell)$, where the measures 
$\tau^{(m)}_k$  are as in $(2.2)$.  
Let $\theta\in (0,1)$ and let positive numbers
 $A=(\log \rho)\|\Omega\|_{L^q(S^{n-1})}\|h\|_{\Delta_s}$, 
$B=\left(1-\beta_m^{-\theta\alpha_m}\right)^{-1}$ be as above.  We define  
$$T_{\rho,m}^*(f)(x) = \sup_{k \in \Bbb Z} \left|\sum_{j=k}^{\infty} 
\tau_j^{(m)}*f(x) \right|.      \tag 4.1  $$   
 Then, for   $p \in (2(1+\theta)/(\theta^2-\theta+2), (1+\theta)/\theta)
 =:I_\theta$    we have  
$$\|T_{\rho,m}^*(f)\|_p \leq  CA
\left(B^{1 +\delta (p)}
+ B^{2/p+1-\theta/2} \right)\|f\|_p, $$  
where   $C$ is independent of $q, s\in (1,2]$,  $\Omega\in 
L^q(S^{n-1})$, $h\in \Delta_s$, $\rho$  and the coefficients of 
the polynomials $P_j$ $(1\leq j\leq d)$. 
\endproclaim  
\demo{Proof}  Let 
$T_\rho^{(m)}(f)=\sum_k\tau^{(m)}_k*f$ be as in Lemma 2.    
Let a function $\varphi$ be as in the definition of $\tau_k^{(m)}$ 
in (2.2).    
 Define $\varphi_k$ by 
$\hat{\varphi}_k(\xi) = 
\varphi\left(\beta_m^k|H_m\pi^d_{s_m}R_m(\xi)|\right)$. 
Let $\delta$ be the delta function as above.  Following  \cite{8}, 
we decompose  
$$\sum_{j=k}^{\infty}\tau^{(m)}_j*f =\varphi_k*T_\rho^{(m)}(f)-
\varphi_k*\left(\sum_{j=-\infty}^{k-1}\tau^{(m)}_j*f \right) + 
  (\delta-\varphi_k)*\left(\sum_{j=k}^{\infty}\tau^{(m)}_j*f\right).$$  
It follows that 
$$ 
T_{\rho,m}^*(f) 
\leq \sup_k\left|\varphi_k*T_\rho^{(m)}(f)\right| 
+ \sum_{j=0}^\infty N_j^{(m)}(f), 
\tag 4.2 
 $$
where   
$N_j^{(m)}(f) = \sup_k\left|\varphi_k*\left(\tau^{(m)}_{k-j-1}*f\right)\right|
+ \sup_k \left|(\delta - \varphi_k)*\left(\tau^{(m)}_{j+k}*f\right)\right|$.  
By Lemma 2 we have     
$$ 
\left\|\sup_k\left|\varphi_k*T_\rho^{(m)}(f)\right|\right\|_p
\leq CAB^{1+\delta(p)} \|f\|_p  
\quad \text{for $p\in (1+\theta, (1+\theta)/\theta)$. }  \tag 4.3 
$$ 
Also, by (2.7) we see that  
$$ 
\|N^{(m)}_j(f)\|_r \leq CAB^{2/r}\|f\|_r 
\quad \text{for $r >1+\theta$. }\tag 4.4
$$ 
\par 
On the other hand, we have 
$$N^{(m)}_j(f) \leq 
\left(\sum_k \left|(\delta - \varphi_k)*\left(\tau^{(m)}_{j+k}*f\right)
\right|^2\right)^{1/2}
+\left(\sum_k \left|\varphi_k*\left(\tau^{(m)}_{k-j-1}*f\right)
\right|^2\right)^{1/2}.$$
Therefore, by the estimates (2.5), (2.6) and  Plancherel's theorem,  
as in \cite{8, p. 820} we see that 
$$\|N^{(m)}_j(f)\|_2 
\leq CA\beta_m^{-\alpha_m j}\left(1-\beta_m^{-2\alpha_m}\right)^{-1/2}\|f\|_2. 
\tag 4.5 $$  
For $p \in I_\theta$ 
we can  find  $r \in (1+\theta, 2(1+\theta)/\theta)$ 
such that $1/p = (1 - \theta)/r + \theta/2$,  so   
 an interpolation between (4.4) and (4.5) implies that
$$\|N^{(m)}_j(f)\|_p \leq CAB^{2(1-\theta)/r}
\left(1-\beta_m^{-2\alpha_m}\right)^{- \theta/2}\beta_m^{-\alpha_m \theta j}
\|f\|_p.
\tag 4.6 $$
Therefore,  by (4.2), (4.3) and (4.6),  for $p \in I_\theta$  we have
$$
\|T_{\rho,m}^*(f)\|_p \leq CA\left(B^{1+\delta(p)} 
+ B^{2(1-\theta)/r +1}
\left(1-\beta_m^{-2\alpha_m}\right)^{- \theta/2} \right)
\|f\|_p.
$$ 
This implies the conclusion of Lemma 6, since 
  $\left(1-\beta_m^{-2\alpha_m}\right)^{-1}\leq B$ and 
 $2(1-\theta)/r +\theta/2+1= 2/p+1-\theta/2$.   
\enddemo
\demo{Proof of Theorem $3$}     
Note that $T^*(f)\leq 2T_{0}^*(f)+2\mu_\rho^*(|f|)$, 
where $T_0^*(f)$  
is defined  by the formula in (4.1) 
with $\{\tau_j^{(m)}\}$ replaced by 
 the sequence $\{\sigma_j\}$ of measures  in (2.1)  
and $\mu_\rho^*=(\mu^{(\ell)})^*$ is as in Lemma 1.   We note that 
$T_0^*(f)\leq \sum_{m=1}^\ell T_{\rho,m}^*(f)$.    
Now,  Lemma 6 implies that  
 $$\|T_{\rho,m}^*(f)\|_p \leq C (\log\rho)
 \left(1-\rho^{-\theta/(2q's')}\right)^{-3}\|\Omega\|_q
 \|h\|_{\Delta_s}\|f\|_p $$
for $p \in I_\theta$.   
 By using this with $\rho=2^{q's'}$, since $\theta\in (0,1)$ is 
 arbitrary,  we can conclude that 
$$\|T^*_{2^{q's'}, m}(f)\|_p \leq C_p(q - 1)^{-1}(s - 1)^{-1}
\|\Omega\|_q\|h\|_{\Delta_s}\|f\|_p $$
for  $p \in (1,\infty)$.  Also, by Lemma 1   
$\mu_\rho^*$  satisfies  a similar estimate when $\rho=2^{q's'}$.  
Collecting results,  we have Theorem 3.  
\enddemo 
\remark{Remark }  
Let       
$$M(f)(x)=\sup_{t>0}t^{-n}
\int_{|y|<t}|f(x-P(y))||\Omega(y')||h(|y|)|\,dy.  $$ 
It is easy to see that $M(f)\leq C\mu_\rho^*(f)$, 
where $C$ is independent of 
$\rho\geq 2$. Therefore, by Lemma 1 we can prove results similar to Theorems 1 
and 2 for the maximal operator $M$.  
In \cite{1},  $L^p$ boundedness of $M$ was proved under the condition 
that $\Omega\in L\log L(S^{n-1})$ and $h\in \Delta_s$ for some $s>1$.  
When $n=d$, $P(y)=y$,
 it is known that $M$ is of weak type $(1,1)$ if $\Omega\in L\log L(S^{n-1})$ 
 and $h$ is bounded (see \cite{5, 4}).  
\endremark

\address
\ 
\newline 
Department of Mathematics 
\newline
Faculty of Education
\newline
Kanazawa University   
\newline
Kanazawa, 920-1192
\newline
Japan
\endaddress

\email
shuichi\@kenroku.kanazawa-u.ac.jp
\endemail 

\enddocument